\input amssym.def 
\input amssym
\magnification=1200
\parindent0pt
\hsize=16 true cm
\baselineskip=13  pt plus .2pt
$ $

\def\C{{\cal C}}

\centerline {\bf Tetrahedral Coxeter groups, large group-actions on 3-manifolds}

\centerline {\bf and equivariant Heegaard splittings}

\bigskip \bigskip

\centerline {Bruno P. Zimmermann}

\medskip

\centerline {Universit\`a degli Studi di Trieste}

\centerline {Dipartimento di Matematica e Geoscienze}

\centerline {34127 Trieste, Italy}

\bigskip \bigskip

{\bf Abstract.}  We consider finite group-actions on closed, orientable and
nonorientable 3-manifolds $M$ which preserve the two handlebodies of a Heegaard
splitting of $M$ of some genus $g > 1$ (maybe interchanging the two handlebodies). 
The maximal possible order of a finite group-action on a handlebody of genus $g>1$
is $12(g-1)$ in the orientation-preserving case and $24(g-1)$ in general, and the
maximal order of a finite group  preserving the Heegaard surface of a Heegaard
splitting of genus $g$ is  $48(g-1)$. This defines a
hierarchy for finite group-actions on 3-manifolds which we discuss in the present
paper; we present various manifolds with an action of type $48(g-1)$ for small
values of $g$, and in particular the unique hyperbolic 3-manifold  with
such an action of smallest possible genus
$g = 6$ (in strong analogy with the Euclidean case of the 3-torus  which has such 
actions for $g = 3$).

\bigskip

{\bf 1. A hierarchy for large finite group-actions on 3-manifolds}

\medskip

The maximal possible order of a finite group $G$ of orientation-preserving
diffeomorphisms of an orientable handlebody of genus $g>1$ is $12(g-1)$ ([Z1],
[MMZ]); for orientation-reversing finite group-actions on an orientable handlebody
and for actions on a nonorientable handlebody, the maximal possible order is
$24(g-1)$; we will always assume $g>1$ in the present paper.

\medskip

Let $G$ be a finite group of diffeomorphisms of a closed, orientable or
nonorientable 3-manifold $M$. We consider Heegaard splittings of $M$  into two
handlebodies of genus $g$ (nonorientable if $M$ is nonorientable)  such that  each
element of $G$ either preserves both handlebodies or interchanges them. The maximal
possible orders are $12(g-1)$, 
$24(g-1)$ or $48(g-1)$, and we distinguish four types of $G$-actions:

\bigskip

{\bf  Definition 1.1}  

\smallskip

1.1.1   {\it The $G$-action is of type $12(g-1)$:} 

\smallskip

$M$ is orientable, $G$ is orientation-preserving and {\it non-interchanging} (i.e.,
preserves both handlebodies of the Heegaard splitting), and $G$ is of maximal
possible order
$12(g-1)$ for such a situation (these actions are called {\it maximally symmetric} 
in various papers, see [Z3,Z5,Z6]).

\medskip

1.1.2  {\it  The $G$-action is of non-interchanging type $24(g-1)$:}

\smallskip

either $M$ is orientable and $G$ orientation-preserving or
$M$ is non-orientable, and $G$ is  of maximal possible order
$24(g-1)$ for such a situation (these actions are called {\it strong maximally
symmetric} in [Z8]).

\medskip

1.1.3  {\it  The $G$-action is of interchanging type $24(g-1)$:}

\smallskip

$M$ is orientable, the subgroup $G_0$ of index 2 preserving both handlebodies
is orientation-preserving, and
$G$ is of  maximal possible order  $24(g-1)$ for such a situation.

\medskip

1.1.4  {\it  The  $G$-action is of type $48(g-1)$:}

\smallskip

$G$ is interchanging, either $M$ is orientable and the subgroup $G_0$ preserving
both handlebodies is orientation-reversing, or $M$ is
non-orientable, and
$G$ is of maximal possible oder $48(g-1)$ for such a situation.

\bigskip

The second largest orders in the four cases are $8(g-1)$, $16(g-1)$ and
$32(g-1)$, then $20(g-1)/3$, $40(g-1)/3$ and
$80(g-1)$, next $6(g-1)$, $12(g-1)$ and
$24(g-1)$ etc.; in the present paper, we consider only the cases of largest
possible orders $12(g-1)$, $24(g-1)$ and
$48(g-1)$.

\medskip

In section 2, we present examples of 3-manifolds for various of these types, 
and in particular the unique hyperbolic 3-manifold of type $48(g-1)$ of
smallest possible genus $g = 6$ . The situation for the orientation-preserving
actions of types 1.1.1 and 1.1.3 is quite flexible and has been considered in
various papers (see [Z3,Z5,Z6]), so in the present paper we concentrate mainly on
the  orientation-reversing actions and actions on non-orientable manifolds of cases
1.1.2 and 1.1.4 where the situation is more rigid. We  finish this section  with a
short discussion of 3-manifolds of type $12(g-1)$, for small values of $g$.

\bigskip

{\bf Theorem 1.2.} ([Z7])  {\sl   The closed orientable  3-manifolds  with a
$G$-action of type $12(g-1)$ and of genus $g = 2$ are  exactly the 3-fold cyclic
branched coverings of the 2-bridge links, the group
$G$ is isomorphic to the dihedral group $\Bbb D_6$ of order 12 and obtained as the
lift of a symmetry group $\Bbb Z_2 \times \Bbb Z_2$ of each such  2-bridge link.}

\bigskip

Examples of such 3-manifolds are the spherical
Poincar\'e homology sphere (the 3-fold branched covering of the torus knot
of type (2,5)), the Euclidean Hantzsche-Wendt manifold (the 3-fold branched
covering of the figure-8 knot, see [Z2]) and the hyperbolic
Matveev-Fomenko-Weeks manifold of smallest volume among all closed hyperbolic
3-manifolds (the 3-fold branched covering of the 2-bridge knot $5_2$).

\medskip

Examples of  3-manifolds of type $12(g-1)$ and genus 3 are the 3-torus and again the
Euclidean Hantzsche-Wendt manifold, of genus 6 the spherical Poincar\'e homology
3-sphere and the hyperbolic Seifert-Weber dodecahedral space, see Corollaries  2.3,
2.4  and 2.5.

\medskip

The finite groups $G$ which admit an action of type $12(g-1)$ and genus $g \le
6$ are $\Bbb D_6$, $\Bbb S_4$, $\Bbb D_3 \times \Bbb D_3$, $\Bbb S_4 \times \Bbb
Z_2$ and $\Bbb A_5$ ([Z6]).

\bigskip

{\bf 2.  Tetrahedral Coxeter groups and large group-actions}

\medskip

{\bf 2.1  Non-interchanging actions of type $24(g-1)$ and actions of type
$48(g-1)$}

\medskip

The following is proved in [Z8].

\bigskip

{\bf Theorem 2.1.}   {\sl  i)  Let $M$ be a closed, irreducible 3-manifold with a
non-interchanging $G$-action of type $24(g-1)$. Then
$M$ is spherical, Euclidean or hyperbolic and obtained as a quotient of the
3-sphere, Euclidean or hyperbolic 3-space by a  normal subgroup $K$ of finite
index, acting freely, in a spherical, Euclidean or hyperbolic Coxeter group 
$C(n,m;2,2;2,3)$ or in a twisted Coxeter group 
$C_\tau(n,m;2,2;3,3)$; the $G$-action is the projection of the Coxeter or twisted
Coxeter group to $M$. Conversely, each such subgroup $K$ gives a $G$-action of
type $24(g-1)$ on the 3-manifold $M = \Bbb H^3/K$.

\smallskip

ii) The $G$-action of type $24(g-1)$ on $M$ extends to a $G$-action of type
$48(g-1)$ if and only $n = m$ and the universal covering group $K$ of $M$ is a
normal subgroup also of the twisted Coxeter group 
$C_\mu(n,n;2,2;2,3)$ or in the doubly-twisted Coxeter group 
$C_{\tau\mu}(n,n;2,2;3,3)$.}

\bigskip

In Theorem 2.1, we use the notation in [Z8] which we now explain.  A {\it 
Coxeter tetrahedron} is a tetrahedron in the 3-sphere, Euclidean or hyperbolic
3-space all of whose dihedral angles are of the form
$\pi/n$ (denoted by a label $n$ of the edge, for some integer $n \ge 2$) and,
moreover, such that at each of the four vertices the three angles of the  adjacent
edges define a spherical triangle (i.e., $1/n_1 + 1/n_2 + 1/n_3 > 1$).  We will
denote such a  Coxeter tetrahedron by
$\C(n,m;a,b;c,d)$ where
$(n,m), (a,b)$ and $(c,d)$ are the labels of pairs of opposite edges, and we denote
by $C(n,m;a,b;c,d)$ the {\it Coxeter group} generated by the reflections in the
four faces of the tetrahedron $\C(n,m;a,b;c,d)$, a properly discontinuous group of
isometries of one of the three geometries. In the following, we list the Coxeter
groups of the various types occuring in Theorem 2.1.

\bigskip

{\bf  2.1.1   The Coxeter groups $C(n,m;2,2;2,3)$:}

\medskip

spherical:  \hskip 5mm $C(2,2;2,2;2,3), \;\;  C(2,3;2,2;2,3), \;\;  C(2,4;2,2;2,3),
\;\; C(2,5;2,2;2,3)$,

\hskip 21mm  $C(3,3;2,2;2,3), \;\;  C(3,4;2,2;2,3), \;\;  C(3,5;2,2;2,3)$;
 
\smallskip

Euclidean:    \hskip 3mm $C(4,4;2,2;2,3)$;

\smallskip

hyperbolic:  \hskip 2mm $C(4,5;2,2;2,3), \;\;   C(5,5;2,2;2,3)$.

\bigskip

A Coxeter tetrahedron  $\C(n,m;2,2;3,3)$ has a rotational symmetry
$\tau$ of order two (an isometric involution) which exchanges the opposite edges
with labels 2 and  3 and inverts the two edges with lables $n$ and $m$. The
involution $\tau$ can be realized by an isometry and defines a {\it twisted Coxeter
group}  $C_\tau(n,m;2,2;3,3)$, a group of isometries containing the Coxeter group 
$C(n,m;2,2;3,3)$ as a subgroup of index two.

\bigskip

{\bf 2.1.2  The twisted Coxeter groups  $C_\tau(n,m;2,2;3,3)$:}

\medskip

spherical:  \hskip 5mm  $C_\tau(2,2;2,2;3,3), \;\;  C_\tau(2,3;2,2;3,3), \;\;
C_\tau(2,4;2,2;3,3)$;

\smallskip

Euclidean:    \hskip 3mm $C_\tau(3,3;2,2;3,3)$;

\smallskip

hyperbolic:  \hskip 2mm $C_\tau(2,5;2,2;3,3), \;\; C_\tau(3,4;2,2;3,3), \;\;
C_\tau(3,5;2,2;3,3)$,  

\hskip  20mm  $C_\tau(4,4;2,2;3,3),  \;\;  C_\tau(4,5;2,2;3,3), \;\; 
C_\tau(5,5;2,2;3,3).$

\bigskip

A Coxeter tetrahedron $\C(n,n;2,2;2,3)$ has a rotational symmetry
$\mu$ of order two (an isometric involution) which exchanges the opposite edges
with labels $n$ and 2 and inverts the two remaining edges with labels 2 and 3. As
before, this defines a {\it twisted Coxeter group} 
$C_\mu(n,n;2,2;2,3)$ containing $C(n,n;2,2;2,3)$ as a subgroup of index 2.

\bigskip

{\bf 2.1.3  The twisted Coxeter groups $C_{\mu}(n,n;2,2;2,3)$:}

\medskip

spherical:  \hskip 5mm  $C_\mu(2,2;2,2;2,3), \;\; C_\mu(3,3;2,2;2,3);$

\smallskip

Euclidean:    \hskip 3mm $C_\mu(4,4;2,2;2,3)$;    

\smallskip

hyperbolic:  \hskip 2mm $C_\mu(5,5;2,2;2,3).$

\bigskip

Finally, a Coxeter tetrahedron $\C(n,n;2,2;3,3)$ has a group $\Bbb Z_2
\times \Bbb Z_2$ of rotational isometries generated by involutions $\tau$ and
$\mu$ as before, and hence defines a {\it doubly-twisted Coxeter group} 
$C_{\tau\mu}(n,n;2,2;3,3)$ containing $C(n,n;2,2;3,3)$ as a subgroup of index 4  
(and both $C_{\tau}(n,n;2,2;3,3)$ and $C_{\mu}(n,n;2,2;3,3)$ as subgroups of index
2).

\bigskip

{\bf 2.1.4  The doubly-twisted  Coxeter groups $C_{\tau\mu}(n,n;2,2;3,3)$:}

\medskip

spherical:  \hskip 5mm  $C_{\tau \mu}(2,2;2,2;3,3)$;

\smallskip

Euclidean:    \hskip 3mm   $C_{\tau \mu}(3,3;2,2;3,3)$;

\smallskip

hyperbolic:  \hskip 2mm  $C_{\tau\mu}(4,4;2,2;3,3), \;\; C_{\tau \mu}(5,5;2,2;3,3)$.

\bigskip

In the following, we will discuss finite-index normal subgroups of small index,
acting freely (i.e., torsion-free in the Euclidan and hyperbolic cases) of various
Coxeter  and tetrahedral groups (their orientation-preserving subgroups). Since 
this requires computational methods, we need presentations of the
various groups.

\bigskip
\vfill  \eject

{\bf 2.2  Presentations of Coxeter and tetrahedral groups}

\medskip

Denoting by  $f_1,f_2,f_3$ and $f_4$ the reflections in the four faces of a
Coxeter polyhedron
$\C(n,m;2,2;c,3)$, the Coxeter group $C(n,m;2,2;c,3)$ has a presentation
$$< f_1,f_2,f_3,f_4 \;\; | \;\;  f_1^2 = f_2^2 = f_3^2 = f_4^2 = 1,$$
$$(f_1f_2)^c = (f_2f_3)^2 = (f_3f_4)^3 = (f_4f_1)^2 =  (f_1f_3)^n = (f_2f_4)^m =
1>.$$ 
A presentation of the twisted
group 
$C_\tau(n,m;2,2;3,3)$ is obtained by adding to this presentation a generator
$\tau$ and the relations
$$\tau^2 = 1, \;\; \tau f_1 \tau^{-1} = f_3, \;\; \tau f_2 \tau^{-1} = f_4.$$   
If $n = m$, for a presentation of $C_\mu(n,n;2,2;c,3)$ one adds a generator $\mu$
and the relations
$$\mu^2 = 1, \;\; \mu f_1 \mu^{-1} = f_2, \;\; \mu f_3 \mu^{-1} = f_4,$$ 
and for a presentation of $C_{\tau \mu}(n,n;2,2;3,3)$ both generators
$\tau$ and $\mu$ with their relations, and also the relation 
$(\tau\mu)^2 = 1$.

\medskip

We consider also the orientation-preserving subgroups of index 2 of the Coxeter
groups. The generators $f_i$ in their presentations denote
rotations now, see [Z5] or [Z6] for such computations of the orbifold fundamental
groups in some of these cases. Representing the 1-skeleton of a tetrahedron by a
square with its two diagonals, a Wirtinger-type representation of the
orbifold-fundamental group is obtained; here the two horizontal edges of the square
have labels
$n$ and $m$, the two vertical edges labels $c$ and $3$ (generators $f_1$ and
$f_4$), the two diagonals labels 2 (generators $f_2$ and $f_3$), and 
similarly for the quotients of the 1-skeleton
of the tetrahedron by the involutions $\tau$ and $\mu$ (represented by rotations
around a vertical and a horizontal axis, so one easily depicts the
singular sets of the quotient orbifolds). In this way on obtains the following
presentations:

\bigskip

the {\it tetrahedral group} $T(n,m;2,2;c,3)$ of index 2 in $C(n,m;2,2;c,3)$:
$$< f_1,f_2,f_3,f_4 \;\; | \;\;  f_1^c = f_2^2 = f_3^2 = f_4^3 = 1,
\;  f_1f_2f_3f_4  =  (f_1f_2)^n = (f_2f_4)^m = 1>;$$  
the {\it twisted tetrahedral group} $T_\tau(n,m;2,2;3,3)$ of index 2 in 
$C_\tau(n,m;2,2;3,3)$: 
$$< f_1,f_2,f_3,f_4 \;\; | \;\;  f_1^2 = f_2^2 = f_3^2 = f_4^3 = 1,
\;  f_1f_2f_3f_4  =  (f_1f_2)^n  = (f_2f_3f_2f_4)^m = 1>;$$ 
the {\it twisted tetrahedral group} $T_\mu(n,n;2,2;c,3)$ of index 2 in 
$C_\mu(n,n;2,2;c,3)$: 
$$< f_1,f_2,f_3,f_4,f_5 \;\; | \;\;  f_1^c = f_2^2 = f_3^2 = f_4^3 = 1, 
\; f_1f_2f_3f_4  = (f_1f_2)^n = 1,$$
$$ f_5^2 = 1,  \; (f_1f_5)^n  =  (f_4f_5)^n = f_3(f_4f_5)f_2(f_4f_5) = 1>;$$ 

the {\it doubly-twisted tetrahedral group}
$T_{\tau \mu}(n,n;2,2;3,3)$ of index 2 in 
$C_{\tau \mu}(n,n;2,2;3,3)$:
$$< f_1,f_2,f_3,f_4,f_5 \;\; | \;\;  f_1^2 = f_2^2 = f_3^2 = f_4^3 = 1,
\; f_1f_2f_3f_4 = (f_1f_2)^n = 1,$$
$$f_5^2 = 1, \; (f_1f_5)^2 = (f_4f_5)^2 = (f_2f_4f_5)^2 = 1>.$$

\medskip

As a typical example, we consider the doubly-twisted
Coxeter group  $C_{\tau \mu}(5,5;2,2;3,3)$ in the next section. 
In the hyperbolic and Euclidean cases, we call a
surjection of a Coxeter group or tetrahedral group {\it admissible} if its kernel
is torsionfree.

\bigskip

{\bf 2.3   Manifolds of type $48(g-1)$ and $24(g-1)$}

\medskip

By Theorem 2.1, we are interested in torsionfree normal subgroups of finite index of
the Coxeter groups 2.1.1 - 2.1.4. Using the presentations in the
previous section, all computations in the following are easily verified by GAP
(which classifies surjections onto  finite groups up to isomorphisms of the image). 
As a typical example, we consider the twisted Coxeter group
$C_\tau(5,5;2,2;3,3)$.

\bigskip

{\bf Theorem 2.2}   {\sl i) There is a unique torsionfree
normal subgroup $K_0$ of smallest possible index 120 in the  twisted Coxeter group 
$C_\tau(5,5;2,2;3,3)$.  The quotient manifold $M_0 = \Bbb
H^3/K_0$  is an orientable hyperbolic 3-manifold of type $48(g-1)$ and genus
$g = 6$, for an action of  $\Bbb S_5  \times \Bbb Z_2$. Since $H_1(M_0) \cong \Bbb
Z^6$, also the ordinary Heegaard genus of $M_0$ is equal to 6.

\smallskip

ii)  The manifold $M_0$ is the unique hyperbolic 3-manifold with an action of type 
$48(g-1)$, and also of non-interchanging type $24(g-1)$, for genera $g \le 6$.

\smallskip

iii) There is a unique torsionfree normal subgroup $K_1$ of smallest possible index
120  in the Coxeter group  $C(5,5;2,2;3,3)$. The quotient $M_1 = \Bbb H^3/K_1$  is
an orientable hyperbolic 3-manifold of type $48(g-1)$ and genus $g = 11$. Since
$H_1(M_1) \cong \Bbb Z^{11}$,  also the ordinary Heegaard genus of $M_1$ is equal
to 11  ($M_1$ is a 2-fold covering of $M_0$)}

\bigskip

{\it Proof.}  i) Since $C_\tau(5,5;2,2;3,3)$ has a finite
subgroup (a vertex group) isomorphic to the extended dodecahedral group 
$\bar \Bbb A_5 \cong \Bbb A_5 \times \Bbb Z_2$ of order 120 (isomorphic to the
extended triangle group $[2,3,5]$ generated by the reflections in the sides of a
hyperbolic triangle with angles $\pi/2, \pi/3$ and $\pi/5$), a torsionfree
subgroup of $C_\tau(5,5;2,2;3,3)$ has index at least 120. Similarly, 
$T_\tau(5,5;2,2;3,3)$ has a vertex group $\Bbb A_5$ and a torsionfree subgroup has
index at least 60.

\medskip

Up to isomorphism of the image (this will always be the convention 
in the following), there are exactly three surjections of  the twisted tetrahedral
group $T_\tau(5,5;2,2;3,3)$ to its vertex group $\Bbb A_5$; the abelianizations
of the kernels of the three surjections  are
$\Bbb Z^6$ and two times  $\Bbb Z_2^4 \times \Bbb Z_4 \times \Bbb Z_5^3$. The three
surjections are admissible (have torsionfree kernel) since, when killing an element
of finite order in a vertex group of the twisted tetrahedron, the twisted
tetrahedral group becomes trivial or of order two. The kernel $K_0$ of the unique
surjection with abelianized kernel $\Bbb Z^6$ is normal also in the twisted Coxeter
group $C_\tau(5,5;2,2;3,3)$, the doubly-twisted tetrahedral group
$T_{\tau\mu}(5,5;2,2;3,3)$, and hence also in the doubly-twisted Coxeter group
$C_{\tau\mu}(5,5;2,2;3,3)$.  By Theorem 2.1, the quotient manifold $M_0 = \Bbb
H^3/K_0$  is a closed orientable hyperbolic 3-manifold of type $48(g-1)$ and genus
$g = 6$.  Since $C_\tau(5,5;2,2;3,3)/K_0  \cong \Bbb A_5 \times \Bbb Z_2$ and there
is no surjection of $C_{\tau\mu}(5,5;2,2;3,3)$ to $\Bbb A_5$, the only remaining
possibility is  $C_{\tau\mu}(5,5;2,2;3,3)/K_0 \cong \Bbb S_5 \times \Bbb Z_2$.

\medskip

There are three surjections of $C_\tau(5,5;2,2;3,3)$ to its vertex group 
$\Bbb A_5 \times \Bbb Z_2$, with abelianizations 
$\Bbb Z^6$,  $\Bbb Z^{12}$ and $\Bbb Z^5 \times \Bbb Z_2^2$, but only the
surjection with kernel is admissible: since the rank of the other two kernels is
larger than 6, they cannot uniformize a 3-manifold with a Heegaard splitting of
genus 6. Hence $K_1$ is the unique torsionfree subgroup of index 120 in 
$C_\tau(5,5;2,2;3,3)$.

\medskip

ii)  By the lists 2.1.3 and 2.1.4 and Theorem 2.1 ii), apart from
$C_{\tau\mu}(5,5;2,2;3,3)$ the other two Coxeter groups to obtain a
hyperbolic 3-manifold of type $48(g-1)$ are $C_{\tau\mu}(4,4;2,2;3,3)$ and
$C_\mu(5,5;2,2;2,3)$.

\medskip

Let $K$ be a torsionfree normal subgroup of $C_\tau(4,4;2,2;3,3)$, with factor
group of order $24(g-1)$. Since the extended octahedral group 
$\Bbb S_4 \times \Bbb Z_2$ of order 48 is a vertex group of $C_\tau(4,4;2,2;3,3)$,
the cases $g = 2, 4$ and 6 are not possible.  Also the case $g = 3$ is not possible
since there is not surjection of 
$C_\tau(4,4;2,2;3,3)$ to  $\Bbb S_4 \times \Bbb Z_2$.

\medskip

Considering $g = 5$, suppose that there is an admissible
surjection of  $C_\tau(4,4;2,2;3,3)$ to a group $G$ of order 96;  
since there are no surjections of  $C(4,4;2,2;3,3)$ onto its vertex group $\Bbb S_4
\times \Bbb Z_2$, its restriction to  $C(4,4;2,2;3,3)$ also surjects onto $G$.
Now $G$ has  $\Bbb S_4 \times \Bbb Z_2$ as a subgroup of index 2; dividing out $\Bbb
Z_2$, it  surjects onto $\Bbb S_4 \times \Bbb Z_2$, and then also
$C(4,4;2,2;3,3)$ surjects onto $\Bbb S_4 \times \Bbb Z_2$. Since no such surjection
exists, this excludes also the case $g = 5$, and hence $g \le 6$ is not possible.

\medskip

Considering the case of $C_\mu(5,5;2,2;2,3)$, there is no surjection 
of $C(5,5;2,2;2,3)$ to its vertex group $\Bbb A_5 \times \Bbb Z_2$, and again 
$g \le 6$ is not possible. 

\medskip

This completes the proof of ii) for the case of actions of type $48(g-1)$;
for the case of actions of non-interchanging type $24(g-1)$ one excludes all other
Coxeter groups in a similar way.

\medskip

iii)  There are exactly three surjections of the tetrahedral group 
$T(5,5;2,2;3,3)$ to $\Bbb A_5$, with abelianized kernels 
$\Bbb Z^{11}$ and two times  $\Bbb Z_2 \times \Bbb Z_3^4 \times \Bbb Z_4^4 \times
\Bbb Z_5^3$.  The kernel $K_1$ of the unique
surjection with abelianized kernel $\Bbb Z^{11}$ is normal also in 
$C(5,5;2,2;3,3)$, $T_\tau(5,5;2,2;3,3)$, $T_\mu(5,5;2,2;3,3)$ and hence also in 
$C_{\tau\mu}(5,5;2,2;3,3)$. By Theorem 2.1 ii), $M_1 = \Bbb H^3/K_1$ is an
orientable 3-manifold to type $48(g-1)$. 

\medskip

Since  there is just one surjection of $C(5,5;2,2;3,3)$ to its vertex group 
$\Bbb A_5 \times \Bbb Z_2$, the kernel $K_1$ is the unique normal subgroup of index
120 in  $C(5,5;2,2;3,3)$.

\medskip

This completes the proof of Theorem 2.2.

\bigskip

The manifold $M_1$ appeared first in [Z4] where 
an explicit geometric description is given.

\bigskip

Next we consider the Coxeter group $C(4,4;2,2;3,3)$. Its subgroup
$T(4,4;2,2;3,3)$ has a unique surjection to ${\rm PSL}(2,7)$, its kernel $K_2$ 
has abelianization  $\Bbb Z^{13}$ and is normal also in $C(4,4;2,2;3,3)$. There are
three surjections of 
$C(4,4;2,2;3,3)$ to ${\rm PSL}(2,7) \times \Bbb Z_2$, exactly one with abelianization 
$\Bbb Z^{13}$, hence its kernel $K_2$ is normal also in the twisted groups
$C_\tau (4,4;2,2;3,3)$, $C_\mu (4,4;2,2;3,3)$ and  
$C_{\tau\mu}(4,4;2,2;3,3)$ and Theorem 2.1 implies:

\bigskip

{\bf Corollary 2.3}   {\sl The quotient manifold  $M_2 = \Bbb H^3/K_2$  is an
orientable  hyperbolic 3-manifold of type $48(g-1)$ and genus $g = 29$.}

\bigskip

As noted in the proof of Theorem 2.2 ii), the Coxeter group $C(5,5;2,2;2,3)$ admits
no surjection onto its vertex group  $\Bbb A_5 \times \Bbb Z_2$. On the other hand,
the tetrahedral group $T(5,5;2,2;2,3)$  has exactly two admissible surjections onto 
its vertex group $\Bbb A_5$, their kernels are conjugate in $C(5,5;2,2;2,3)$,
normal in $T_\mu(5,5;2,2;2,3)$ with factor group
$\Bbb S_5$, and uniformize the Seifert-Weber hyperbolic
dodecahedral 3-manifold ([WS]).  By [M], $\Bbb S_5$ is in fact  the full
isometry group of the Seifert-Weber manifold which has no orientation-reversing
isometries.  

\bigskip

{\bf Corollary 2.4}   {\sl The Seifert-Weber hyperbolic dodecahedral 3-manifold is a
closed orientable 3-manifold of interchanging type $24(g-1)$ and genus $g = 6$, for
the action of its isometry group $\Bbb S_5$.}

\bigskip

There are three surjections of the tetrahedral group $T(5,5;2,2;2,3)$
to ${\rm PSL}(2,19)$, all admissible, and exactly one kernel 
$K_3$ has infinite abelianization $\Bbb Z^{56}$
and is normal also in the Coxeter group $C(5,5;2,2;2,3)$. There is exactly
one surjection of $C(5,5;2,2;2,3)$ to ${\rm PSL}(2,11) \times \Bbb Z_2$, with kernel $K_3$, 
hence $K_3$ is normal also in
the twisted Coxeter group $C_\mu(5,5;2,2;2,3)$ and Theorem 2.2 implies:

\bigskip

{\bf Corollary 2.5}   {\sl The quotient manifold  $M_3 = \Bbb H^3/K_3$  is an
orientable  hyperbolic 3-manifold of type $48(g-1)$ and genus $g = 286$.}

\bigskip

Infinite series of finite quotients of the hyperbolic Coxeter group
$C(5,5;2,2;2,3)$ are considered in the papers [JL] and [P].

\medskip

It is shown in [JM] that the twisted Coxeter group $C_\tau(5,2;2,2;3,3)$
has a  unique  torsionfree normal subgroup of smallest possible index  2640,
with factor group  ${\rm PGL}(2,11)
\times \Bbb Z_2$, which uniformizes an orientable hyperbolic 3-manifold of
non-interchanging type $24(g-1)$ and genus $g = 111$, for an action of ${\rm
PGL}(2,11) \times \Bbb Z_2$.

\bigskip

Next we discuss the  case of the Euclidean Coxeter group
$C_\tau(3,3;2,2;3,3)$  (with strong analogies with the hyperbolic case 
of $C_\tau(5,5;2,2;3,3)$ in Theorem 2.2).

\medskip

There are five surjections of the tetrahedral group
$T(3,3;2,2;3,3)$ to its vertex group $\Bbb A_4$, the abelianizations of the kernels
are  $\Bbb Z^3$, two times $\Bbb Z_4^2$ and two times $\Bbb Z_2^5$. The kernel with 
abelianization $\Bbb Z^3$ uniformizes the
3-torus and is normal also in $C(3,3;2,2;3,3)$, $T_\tau(3,3;2,2;3,3)$ and hence
$C_{\tau\mu}(3,3;2,2;3,3)$.  The two surjections with abelianized kernel  
$\Bbb Z_4^2$ are conjugate in 
$C(3,3;2,2;3,3)$ and uniformize the Hantzsche-Wendt manifold (the only of
the six orientable Euclidean 3-manifolds with homology $\Bbb Z_4^2$, see
[W]), and the remaining two surjections are not admissible.

\bigskip

{\bf Corollary 2.6}   {\sl i) The 3-torus is a 3-manifold of type $48(g-1)$ and
genus $g = 3$ which is also its ordinary Heegaard
genus.

\smallskip

ii) The Euclidean Hantzsche-Wendt manifold is of interchanging type $24(g-1)$ and
genus $g = 3$, for the action of its orientation-preserving
isometry group $\Bbb S_4 \times \Bbb Z_2$ (by section 1, it is also a 3-manifold of
type $12(g-1)$ and genus $g = 2$, for an action of the dihedral group $\Bbb D_6$ of
order 12).}

\bigskip

By [Z2], the full isometry group of the Hantzsche-Wendt manifold has order 96
but the orientation-reversing elements do not preserve the Heegaard splitting of
genus 3 of Theorem 2.4.

\medskip

The 3-torus is a manifold of type $48(g-1)$ and genus $g =
3$ in still a different way.  The Euclidean tetrahedral group $T(4,4;2,2;2,3)$ has
three surjections to its vertex group $\Bbb S_4$, with abelianized kernels $\Bbb
Z^3$ and two times $\Bbb Z_2^2 \times \Bbb Z_4$. The kernel with abelianization
$\Bbb Z^3$ is normal also in the Coxeter group $C(4,4;2,2;2,3)$, the twisted
tetrahedral group 
$T_\mu(4,4;2,2;2,3)$ and hence in the twisted Coxeter group
$C_\mu(4,4;2,2;2,3)$, it uniformizes the 3-torus which is again a 3-manifold 
of type $48(g-1)$ (but for an action not equivalent to
the action arising from $T(3,3;2,2;3,3)$).

\bigskip

As a spherical case, we consider the Coxeter group  $C(5,3;2,2;2,3)$.
There is no surjection of $C(5,3;2,2;2,3)$ to its vertex group 
$\Bbb A_5 \times \Bbb Z_2$, and there are two surjections of the
tetrahedral group $T(5,3;2,2;2,3)$ to its vertex group $\Bbb A_5$; both kernels have
trivial abelianization, are conjugate in  $C(5,3;2,2;2,3)$ and  uniformize the
spherical Poincar\'e homology 3-sphere.

\bigskip

{\bf Corollary 2.7}   {\sl The spherical Poincar\'e homology 3-sphere is a
3-manifold of type $12(g-1)$ and genus $g = 6$, for an action of $\Bbb A_5$.}

\bigskip

The case of non-orientable manifolds is more
elusive. It is shown in [CMT] that, for sufficiently large $n$, the alternating
group $\Bbb A_n$ is a quotient of $C_\tau(5,2;2,2;3,3)$ by a torsionfree normal
subgroup, and such a subgroup uniformizes a non-orientable 3-manifold: since $\Bbb
A_n$ is simple, the orientation-preserving subgroup  $T_\tau(5,2;2,2;3,3)$ of index
2 surjects onto $\Bbb A_n$, and hence the kernel contains an orientation-reversing
element.

\bigskip

{\bf Corollary 2.8}   {\sl For all sufficiently large $n$, there is a 
non-orientable hyperbolic 3-manifold of non-interchanging type $24(g-1)$, for an 
action of the alternating group $\Bbb A_n$.}

\bigskip

Some other finite simple quotients of the nine hyperbolic Coxeter groups  are listed
in [H], and for the Coxeter groups of type $C(n,m;2,2;2,3)$ these define
non-orientable hyperbolic 3-manifolds of non-interchanging type
$24(g-1)$.  At present, we don't know explicit examples of small genus of
non-orientable manifolds of non-interchanging type $24(g-1)$, and no example  
of type $48(g-1)$.

\bigskip \bigskip

\centerline {\bf References}

\bigskip

\item {[CMT]}  M. Conder, G. Martin, A. Torstensson, {\it  Maximal symmetry groups
of hyperbolic 3-manifolds.}  New Zealand J. Math. 35 (2006), 37-62

\smallskip

\item {[H]}  S.P. Humphries, {\it  Quotients of Coxeter complexes, fundamental
groupoids and regular graphs.}  Math. Z.  (1994), 247-273

\smallskip

\item {[JL]} G.A. Jones, C.D. Long,  {\it Epimorphic images of the [5,3,5] Coxeter
group.}  Math. Z. 275 (2013), 167-183

\smallskip

\item {[JM]} G.A. Jones, A.D. Mednykh,  {\it Three-dimensional hyperbolic manifolds
with large isometry groups.}  Preprint 2003

\smallskip

\item {[M]} A.D. Mednykh,  {\it On the isometry group of the Seifert-Weber
dodecahedron.}  Siberian Math. J. 28 (1987), 798-806

\smallskip

\item {[MMZ]} D. McCullough, A. Miller, B. Zimmermann,  {\it Group actions on
handlebodies.}  Proc. London Math. Soc.  59   (1989), 373-415

\smallskip

\item {[P]} L. Paoluzzi,  {\it $PSL(2,q)$ quotients of some hyperbolic
tetrahedral and Coxeter groups.}  Comm. Algebra  26 (1998), 759-778

\smallskip

\item {[WS]} C. Weber, H. Seifert,  {\it  Die beiden Dodekaederr\"aume.}
 Math.Z. 37,  237-253  (1933)

\smallskip

\item {[W]}  J. Wolf, {\it Spaces of Constant Curvature.} Publish or Perish, Boston
1974

\smallskip

\item {[Z1]} B. Zimmermann,  {\it \"Uber Abbildungsklassen von Henkelk\"orpern.} 
Arch. Math. 33  (1979),  379-382

\smallskip

\item {[Z2]} B. Zimmermann,  {\it On  the Hantzsche-Wendt manifold.} Monatsh.
Math. 110  (1990), 321-327 

\smallskip

\item {[Z3]}  B. Zimmermann, {\it Finite group actions on handlebodies and
equivariant Heegaard genus for 3-manifold.}  Top. Appl. 43 (1992), 263-274

\smallskip

\item {[Z4]}  B. Zimmermann, {\it On a hyperbolic 3-manifold with some special
properties.} Math. Proc. Camb. Phil. Soc. 113 (1993), 87-90

\smallskip

\item {[Z5]}  B. Zimmermann, {\it Hurwitz groups and finite group actions on
hyperbolic 3-mnifolds.}  J. London Math. Soc. 52 (1995), 199-208

\smallskip

\item {[Z6]} B. Zimmermann,  {\it Genus actions of finite groups on
3-manifolds.}  Mich. Math. J. 43 (1996),  593-610

\smallskip

\item {[Z7]}  B. Zimmermann, {\it Determining knots and links by cyclic branched
coverings.}  Geom. Ded. 66  (1997),  149-157

\smallskip

\item {[Z8]} B. Zimmermann,  {\it On large orientation-reversing group-actions on
3-manifolds and equivariant Heegaard decompositions.}  Monatsh. Math. 191 (2020), 
437-444

\bye